\documentclass[12pt]{amsart}
\usepackage[english]{babel}
\usepackage[latin1]{inputenc}
\usepackage{hyperref}
 \usepackage{stmaryrd}
\usepackage{amssymb,amsfonts,amsmath,amsthm}
\usepackage{graphicx}
\usepackage[all]{xy}
\usepackage{enumerate}

\usepackage{xcolor}

\usepackage{mathrsfs} 

\newcommand{\ol}[1]{\overline{#1}}

\newcommand{\M}{\mathrm{M}}

\newcommand{\End}{\mathrm{End}}

\newcommand{\C}{\mathbb{C}}

\newtheorem*{prop*}{Proposition}

\theoremstyle{definition}

\newtheorem*{thm*}{Theorem}
\theoremstyle{remark}

\begin{document}
	\title[Fixed point branes, singular loci and mirror symmetry]
	{Fixed point branes, singular loci \\ and mirror symmetry}
	\author[Ana Pe{\'o}n-Nieto]{Ana Pe{\'o}n-Nieto}
	\address{Universit\'e de Gen\`eve\\
		Section de Math\'ematiques\\ Villa Battelle\\
		Route de Drize 7\\ 1227 Carouge\\ Switzerland}
	\email{ana.peon-nieto@unige.ch}
	
	\date{04/07/2019}
	\begin{abstract}
	This is the extended abstract of the talk given at the workshop ``Geometry and physics of Higgs bundles'', held at the Mathematisches Forschungsinstitut Oberwolfach in May 2019. 	 
	\end{abstract}
	\maketitle
	%
	%

	Higgs bundles \cite{SDE} are a vast subject with many ramifications, of which we hereby focus on mirror symmetry. More specifically, we focus on two pieces of work by the author and collaborators \cite{FP,FGOP1}, where  we consider related families of branes  meaningful in topologial mirror symmetry and towards the  understanding of
	mirror symmetry beyond the generic locus.
	
	Let $\M_X$ be the moduli space of Higgs bundles on a Riemann surface $X$, that is, the scheme parametrizing pairs $(E,\phi)$ where $E\longrightarrow X$ is a holomorphic vector bundle of rank $n$ and $\phi\in H^0(X,\End(E)\otimes K_X)$, with $K_X$ the canonical bundle of $X$ \cite{SDE,Ni}. The non abelian Hodge correspondence \cite{Simp} endows it with 3 complex structures $I$ (naturally induced by the complex structure of $X$), $J$ (coming from the character variety) and $K:=I\circ J$, which underlie a hyperk\"ahler structure.
	
	A key tool in the study of  $\M_X$ is the Hitchin map, associating to a pair $(E,\phi)$ the characteristic polynomial of $\phi$. The fibers of this map can be identified with the Jacobian of a suitable spectral cover of $X$ \cite{Duke, DG}, and are  Lagrangian for the $I$-homolorphic  symplectic form $\Omega_I$. So after a hyperk\"ahler rotation, they become special Lagrangian. Thus we have a hyperk\"ahler (hence Calabi--Yau) manifold admitting a special Lagrangian torus fibration, which together with autoduality of Jacobian varieties gives an example of SYZ-mirror symmetry (between $\M_X$ and $\M_X$).	More generally, for any  complex reductive Lie group $G$, the moduli space of $G$-Higgs bundles $\M_G$ is mirror to $\M_{\check{G}}$, where $\check{G}$ denotes the Langlands dual group of $G$  \cite{DP,HT}. Indeed, their Hitchin systems are dual Lagrangian torus fibrations, duality realised by  a Fourier--Mukai equivalence $D(\M_G)\cong D(\M_{\check{G}})$ \cite{DP}.

	In terms of the homological mirror symmetry  conjecture \cite{K}, predictions from physics allow us to interpret mirror symmetry for Hitchin systems as a correspondence between branes \cite{KW}. These can be of type $A$ (that is  Lagrangian submanifolds  with a flat bundle)  or $B$  (complex submanifolds with a holomorphic bundle) in each of the distinguished K\"ahler  structures. Particularly important are $BBB$-branes, and their dual $BAA$-branes.

	All of the above applies over smooth spectral curves. But the global statement of mirror symmetry for Hitchin systems is not well understood. In \cite{FP}, we explore this phenomenon over the singular locus of $\M_X$. From another point of view, a way to evidence for global dualities is by producing global invariants \cite{HMP}. A remarkable example of the above is the equality of the (stringy) E-polynomials of $M_{\mathrm{PGL}(n,\C)}$ and   $M_{\mathrm{SL}(n,\C)}$  \cite{HT,GWZ}. The former can be expressed in terms E-polynomials of some branes in $M_{\mathrm{SL}(n,\C)}$,  investigated in \cite{FGOP1}, whose geometry brings us back to the branes of singular loci from \cite{FP}.

	%
	%
	\section*{Singular loci and their conjectural duals}
	%
	%
	The singular locus of $\M_X$ may be covered with $BBB$-branes  \cite{FP}. These are given by Higgs bundles whose structure group reduces to the Levi subgroup $L_{\ol{n}}$ of the standard parabolic $P_{\ol{n}}$ associated to the partition  $\ol{n}$.  Hypercomplexity of these subspaces  $\M_X(L_{\ol{r}})\subset\M_X$  follows from the non abelian Hodge correspondence. For the partition $(1,\dots,1)$ we produce a family of hyperholomophic bundles $\mathscr{L}$  on $\M_X(L_{\ol{1}})$ paramerized by $\mathrm{Jac}(X)$.

	We next investigate what the duals of the branes $(M(L_{\ol{1}}),\mathscr{L})$ should be, the idea being the potential existence of a stacky Fourier--Mukai transform of $\mathscr{L}$. 
	
	With this in mind, given $F_1,\dots, F_s$ stable vector bundles of rank $\mathrm{rk}(F_i)=n_i$ (where $\sum_in_i=n$), we consider $\M_X(P_{\ol{n}})\subset\M_X$ given by Higgs bundles whose structure group reduces to the parabolic subgroup $P_{\ol{n}}$. Letting $U_{\ol{n}}<P_{\ol{n}}$ be the unipotent radical, let 
	\begin{equation}\label{eq:BAAsingular}
	D_{\ol{n}}^X(\ol{F})=\left\{(E,\phi)\in\M_X(P_{\ol{n}})\,:\, E/U\cong\bigoplus_i F_i\right\}.
	\end{equation}
	Under suitable conditions on the $F_i$'s, $D_{\ol{n}}^X(\ol{F})$ is Lagrangian.  In relation with singular loci, these Lagrangians look very much like  Fourier--Mukai transforms. This, together with the existence of such a transform on the level of generic points of both branes, led us to conjecture the duality between them.
	
	\section*{Branes in topological mirror symmetry}
	Let $\Gamma\subset\mathrm{Jac}(X)$ be the $n$-torsion subgroup. Given $\gamma\in\Gamma$, we consider  its action on $\M_X$ by tensorization $\gamma\cdot (E,\phi)=(E\otimes\gamma, \phi)$. The fixed point set $\M^\gamma$ is easily seen to be hypercomplex \cite{NR}, thus defining the support of a brane of type $BBB$. These are the branes appearing in the expression of the (stringy) E-polynomial of $\M_{\mathrm{PGL}(n,\C)}$, which motivated us to further study them.
	
	When $\gamma\in \Gamma$ is of maximal order, we produce  hyperholomorphic bundles $\mathscr{F}$ on $\M^\gamma$ parametrised by $F\in\mathrm{Jac}(X_\gamma)$, where  $p_\gamma\,:\,X_\gamma\to X$ is the \'etale cover associated with $\gamma$. Since generic spectral curves are generically integral (although always singular), we may generically perform a Fourier--Mukai transform to compute the support of the dual branes $D_\gamma^X({F})$ using the constructions in \cite{A}. 
	Unfortunately, we do not have a good global understanding of $D^X_\gamma(F)$. In fact, in order to prove they are Lagrangian, we check that pullback by $p_\gamma$ yields a local isomorphism from $D_\gamma^X(F)$ into $D_{(1,\dots, 1)}^{X_{\gamma}}(\ol{F})$ (see \eqref{eq:BAAsingular}). 
	
	More generally, given $\gamma\in \Gamma$ of order $m|n$, and a $F$ a rank $n':=n/m$ vector bundle on $m$, we define an isotropic  subscheme $
	D_{\gamma}^{X}({F})
	$ 	which is Lagrangian under suitable conditions on $F$. Again, isotropicity follows from the fact that 
	$$p_\gamma^*\,:\, D_\gamma(F)^{X} \longrightarrow D_{(n',\dots,n')}^{X_{\gamma}}(\ol{F})
	$$ 
	is a local isomorphism. Thus,  in a certain sense, the geometry of singular loci and their conjectural duals determines that of fixed point branes and their duals. 
	
	%
	%
	\section*{Future directions}
	
	An important gap in our constructions, both in \cite{FP} and \cite{FGOP1}, is the lack of a construction of hyperholomorphic bundles beyond the simplest case (Higgs bundles for Cartan subgroups and points fixed by maximal order line bundles).  An immediate improvement (in progress) would consist in extending these results. 
	
	Another natural step is to compute stacky Fourier--Mukai duals of singular loci and check whether stability is preserved. We will undertake this in the rank two case \cite{FP2}, using \cite{M}. 
	
	Finally, it should be pointed out that the study of $\M^\gamma$ is but a first step in a long term project aiming at understanding topological mirror symmetry in terms of branes. To this aim, we carefully study in \cite{FGOP2} the rank two case, and verify that a variation of the hyperholomorphic bundle allows to recover all the meaningful points in the nilpotent cone. 
	
\end{document}